\documentclass[12pt]{article}

\usepackage[english]{babel}%
\usepackage{graphicx}
\usepackage{amsfonts}

\begin{document}

\newtheorem{theorem}{Theorem}
\newtheorem{cor}{Corrolary}
\newtheorem{example}{Example}

\renewcommand{\thesection}{\arabic{section}}
\renewcommand{\thetheorem}{\arabic{section}.\arabic{theorem}}
\renewcommand{\theequation}{\arabic{section}.\arabic{equation}}
\renewcommand{\theexample}{\arabic{section}.\arabic{example}}
\renewcommand{\thecor}{\arabic{section}.\arabic{cor}}

\def\qed{{\bf \hfill $\Box$}\endtrivlist}
\def\endsolve{{\bf \hfill $\bigtriangleup$}\endtrivlist}

\begin{center}{\Large Length of a Full Steiner Tree as a Function of Terminal Coordinates} \end{center}

\begin{center}{\large Alexei Yu. Uteshev, Elizaveta A. Semenova\footnote{\{alexeiuteshev,semenova.elissaveta\}@gmail.com}}  \end{center}

\begin{center}{\it St. Petersburg State University, St. Petersburg, Russia} \end{center}

%\hline

Given the coordinates of the terminals $ \{(x_j,y_j)\}_{j=1}^n $ of the full Euclidean Steiner tree, its length equals
$$ \left| \sum_{j=1}^n z_j U_j \right| \, , $$
where $ \{z_j:=x_j+  \mathbf i y_j\}_{j=1}^n $  and $ \{U_j\}_{j=1}^n $ are suitably chosen $ 6 $th roots of unity. We also extend  this result for the cost of the optimal
Weber networks which are topologically equivalent to some full Steiner trees.

\section{Introduction}
\setcounter{equation}{0}
\setcounter{theorem}{0}
\setcounter{example}{0}

The problem of construction of the shortest possible network interconnecting the set of points (further referred to as \textbf{terminals}) $ \mathbb P=\{P_j\}_{j=1}^n \subset \mathbb R^d, n \ge 3,d\ge 2 $   is known as the \textbf{Euclidean}  \textbf{Steiner minimal tree problem} (\textbf{SMT} problem). We will treat here only the planar version of the problem, i.e. $ d=2 $. For this case, some specializations of the set $ \mathbb P $ guarantee the solution to the problem in the form of a minimum spanning tree. However, for other configurations of $ \mathbb P $,  the optimal network contains some extra vertices (junctions) known as the Fermat-Torricelli or \textbf{Steiner points}. Every such a point in the optimal network possesses the degree $ 3 $, and has the adjacent edges meeting at $ 2\pi/3 $.

The SMT problem has a long and intriguing history \cite{Brazil_2014}. It is known to be NP-hard. Establishing the topology of the minimal tree or even estimation of its length for the general set $ \mathbb P $ is the subject of numerous and inventive research efforts. In the variety of potential solutions, it is possible to distinguish a type of subtrees with a relatively ordered structure.
A Steiner tree with $ n $ terminals is called a \textbf{full Steiner tree} if it contains $ n-2 $ Steiner points. In this case, each terminal has degree equal to $ 1 $.
 Every Steiner minimal tree can be (uniquely) represented as a union of full Steiner trees constructed for the subsets of the set $ \mathbb P $. Any pair of these subtrees might have at most one common terminal.

The first aim of the present paper is express the length of a full Steiner tree as an explicit function of the terminal coordinates. In \cite{Booth} this length is computed  via the lengths of the segments $ |P_jP_{j+1}| $ (for the appropriate numeration of terminals).
In \cite{Weng} the expression for the length was obtained in the so-called \emph{hexagonal coordinates of terminals} which are claimed to be \emph{more natural than the common Cartesian coordinate system}. We intend to dispute this claim by presenting the formula for the length in terms of Cartesian coordinates. For the aim of deducing the formula from the Abstract, we will place the problem into the complex plane. In Section \ref{SLength} we present the proof of the length formula and discuss its relation to one Maxwell's result \cite{Maxwell_1864}.

A natural generalization of the SMT problem is
the (multifacility) \textbf{Weber problem}. It is stated as that of location of the  given number  $ \ell \ge 2 $ of points (named \textbf{facilities}) $ \{W_{r}\}_{r=1}^{\ell} \subset \mathbb R^d $ connected to the terminals of the set $ \mathbb P $ that solve the optimization problem
\begin{eqnarray}
\min_{\{W_{1},\dots, W_{\ell}\}\subset \mathbb R^d} \left\{ \sum_{j=1}^n \sum_{r=1}^{\ell} m_{rj} |W_rP_j| +
\sum_{k=1}^{\ell} \sum_{r=k+1}^{\ell-1} \widetilde m_{rk} |W_rW_k|
\right\} \, ;
\label{F_Weber_m}
\end{eqnarray}
here some of the \textbf{weights} $ m_{ij} $ and $ \widetilde m_{ik} $ might be zero.
The value (\ref{F_Weber_m}) will be referred to as the \textbf{minimal cost of the network}. It can be expected that this problem should be more complicated than the SMT one. It really is even for the planar case, i.e. for $ d=2 $. In contrast with the SMT, the optimal Weber network might contain facilities of degree greater than $ 3 $. However, we conjecture that in a restricted version of the planar Weber problem, one may expect the existence
of the counterpart for the Steiner tree length formula expressing the minimal cost by radicals with respect to the problem parameters of
(terminal coordinates and weights).  This relates to the networks topologically equivalent to some full Steiner trees, i.e. the networks where each terminal  has degree equal to $1 $ while the degree of each facility equal to $ 3 $.
We discuss this issue in Section \ref{SWeber_L}.

{\bf Notation.} We set  $ z_j:=x_j+ \mathbf i y_j $ and denote the $ 6 $th roots of unity
\begin{equation}
\upsilon_k:=\cos \frac{\pi k}{3} + \mathbf i \sin \frac{\pi k}{3} \quad \mbox{for} \ k\in \{0,1,\dots,5\} \, ,
 \label{Roots_6}
\end{equation}
while the third roots of unity as
\begin{equation}
\varepsilon_0=1, \varepsilon_1=-\frac{1}{2}+ \mathbf i \frac{\sqrt{3}}{2}, \ \varepsilon_2=-\frac{1}{2}- \mathbf i \frac{\sqrt{3}}{2} ;
 \label{Roots_3}
\end{equation}
Evidently, $ \upsilon_0=1=-\upsilon_3, \upsilon_1=-\upsilon_4=-\varepsilon_2,\upsilon_2=-\upsilon_5=\varepsilon_1 $.

$ \mathfrak L $ and $ \mathfrak C $ stand for respectively minimal length of the Steiner network and minimal cost of the Weber network.

For the convenience of the reader, we place here the Classical Geometry result which we often refer to in what follows:

\begin{theorem}[Ptolemy] \label{ThPtolemy}
A quadrilateral is inscribable in a circle if and only if the product of the lengths of its diagonals is equal to the sum of the products of the lengths of the pairs of opposite sides.
\end{theorem}

\section{Hints}\label{SHints}
\setcounter{equation}{0}
\setcounter{theorem}{0}
\setcounter{example}{0}

%Some related computations can be found in \cite{Uteshev_Steiner}.

Geometric conditions for the existence of a full Steiner tree for the case $ n=3 $ and $ n=4 $ terminals are well known,
  and we assume them to be fulfilled for the terminal sets treated in the present section. In \cite{Uteshev_Steiner} these conditions are converted into
  analytical form, aside with formulas for Steiner point coordinates and the tree length. This provides one with an opportunity to verify directly the test examples of this and the next sections.

For the case of $ 3 $ terminals, the Torricelli-Simpson construction of the Steiner (Fermat-Torricelli) point
is based on finding an extra point $ Q_1 $ in the plane yielding the equilateral triangle $ P_1P_2Q_1 $ (Fig.1).

\vskip5mm

\begin{minipage}[t]{90mm}
\begin{center}
\includegraphics[width=80mm]{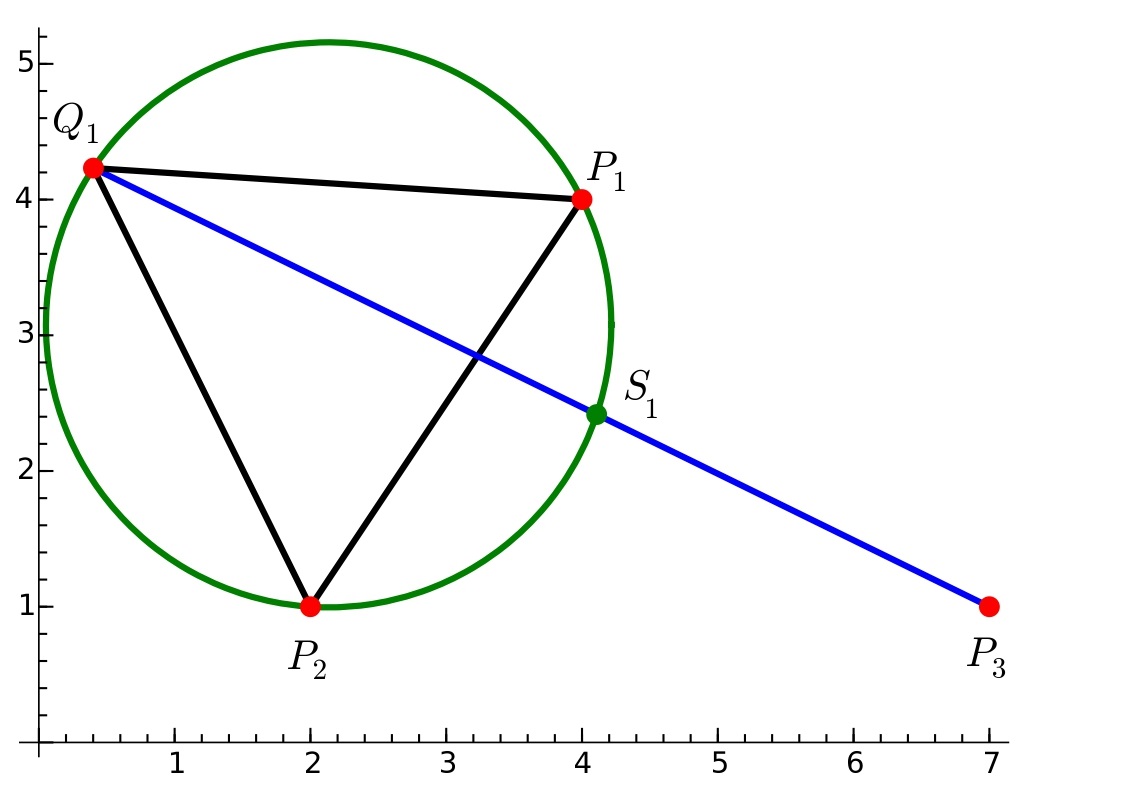}

(a)
\end{center}
\end{minipage}
\hfill
\begin{minipage}[t]{90mm}
\begin{center}
\includegraphics[width=80mm]{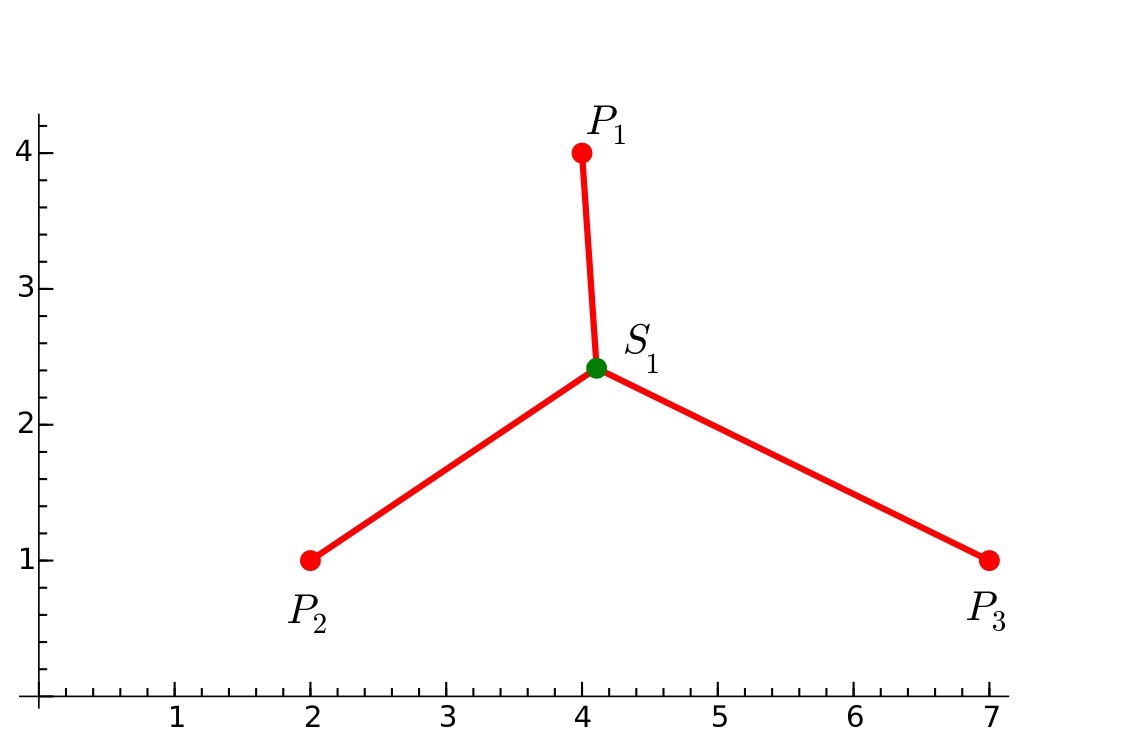}

   (b)
\end{center}
\end{minipage}

\begin{center}
Figure~1. Steiner tree construction for three terminals
\end{center}

\vskip5mm

The point of intersection of the segment $ Q_1P_3 $ with the circle circumscribing this triangle  is the Steiner point $ S_1 $. Due to Theorem \ref{ThPtolemy}, one gets
$$ |Q_1S_1|= |P_1S_1|+|P_2S_1| \, ,$$
and therefore the length of the SMT connecting the terminals $ P_1, P_2, P_3 $ equals $ |Q_1P_3| $. The coordinates of the point $ Q_1 $ can be easily determined. There are two possibilities for its location (on each side of the line $ P_1P_2 $), and we need the variant opposite to $ P_3 $. For a configuration similar to that displayed in Fig.~1 (i.e. the triangle vertices are numbered counterclockwise),
one gets:
\begin{equation}
Q_1:=(q_{1x},q_{1y})=\left(\frac{1}{2}x_1+\frac{1}{2}x_2-\frac{\sqrt{3}}{2}y_1+\frac{\sqrt{3}}{2}y_2\, ,\ \frac{\sqrt{3}}{2}x_1-\frac{\sqrt{3}}{2}x_2+\frac{1}{2}y_1+\frac{1}{2}y_2 \right) \, .
\label{Q1}
\end{equation}
or, alternatively, in the complex plane
$$
\mathfrak q_1:=q_{1x}+\mathbf i q_{1y}=\left(\frac{1}{2}+ \mathbf i \frac{\sqrt{3}}{2}\right)x_1+\left(\frac{1}{2}- \mathbf i \frac{\sqrt{3}}{2}\right)x_2+
\left(  -\frac{\sqrt{3}}{2}+\mathbf i \frac{1}{2} \right)y_1+ \left(\frac{\sqrt{3}}{2}+\mathbf i \frac{1}{2} \right)y_2
$$
$$
=\left(\frac{1}{2}+ \mathbf i \frac{\sqrt{3}}{2}\right)z_1+\left(\frac{1}{2}- \mathbf i \frac{\sqrt{3}}{2}\right)z_2
=\upsilon_1 z_1+\upsilon_5z_2=-\varepsilon_2z_1-
\varepsilon_1z_2 \, .
$$
For this case, the length of the (full) Steiner tree equals
\begin{equation}
 \mathfrak L = |Q_1P_3 |=|z_3- \mathfrak q_1|=|\upsilon_4z_1+\upsilon_2z_2+z_3|=|z_3+\varepsilon_2 z_1 + \varepsilon_1 z_2 | \, .
 \label{Length3}
\end{equation}

\begin{example} \label{Ex0}  For the terminals
$$
P_1=(4,4), \ P_2= (2,1), P_3=(7,1) \, ,
$$
(Fig.~1),
formula (\ref{Length3}) yields
$$ \mathfrak L =|7+\mathbf i + \varepsilon_2 (4+4 \mathbf i) + \varepsilon_1 (2+ \mathbf i) | = | 4+3\sqrt{3}/2-\mathbf i(3/2 + \sqrt{3})  |=\sqrt{28+15\sqrt{3}} \approx 7.347160\, .$$
\end{example}

For the case of $ 4 $ terminals $ \{P_j\}_{j=1}^4 $, the geometrical algorithm for Steiner tree construction \cite{Gilbert&Pollak} with two Steiner points
in the topology $ \begin{array}{c} P_1 \\ P_2 \end{array} S_1 S_2 \begin{array}{c} P_4 \\ P_3 \end{array} $, is illuminated in Fig.~2.
\begin{figure*}
\begin{center}
\includegraphics[width=0.60\textwidth]{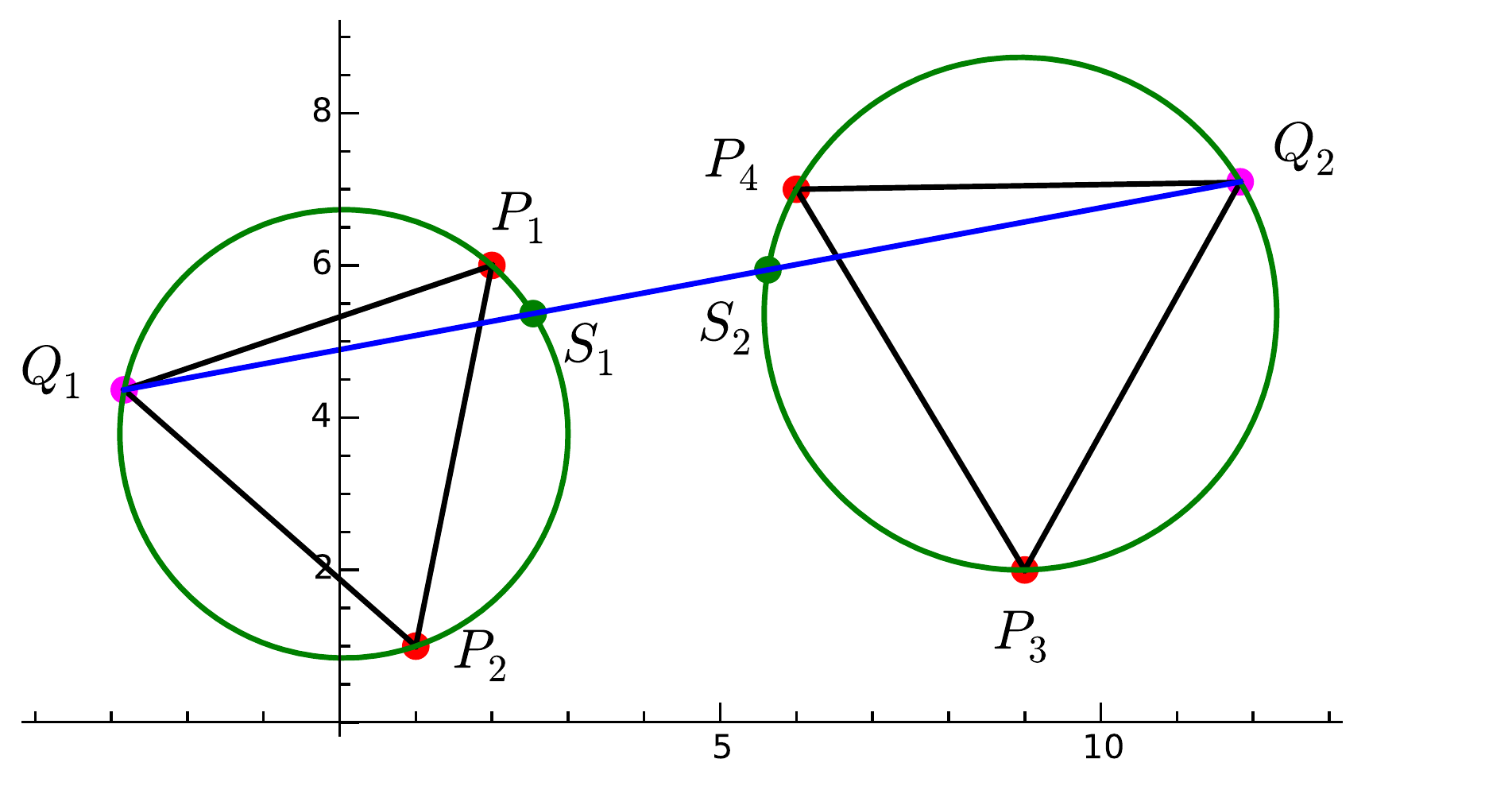}
%\caption{}
\end{center}
\begin{center}
\includegraphics[width=0.60\textwidth]{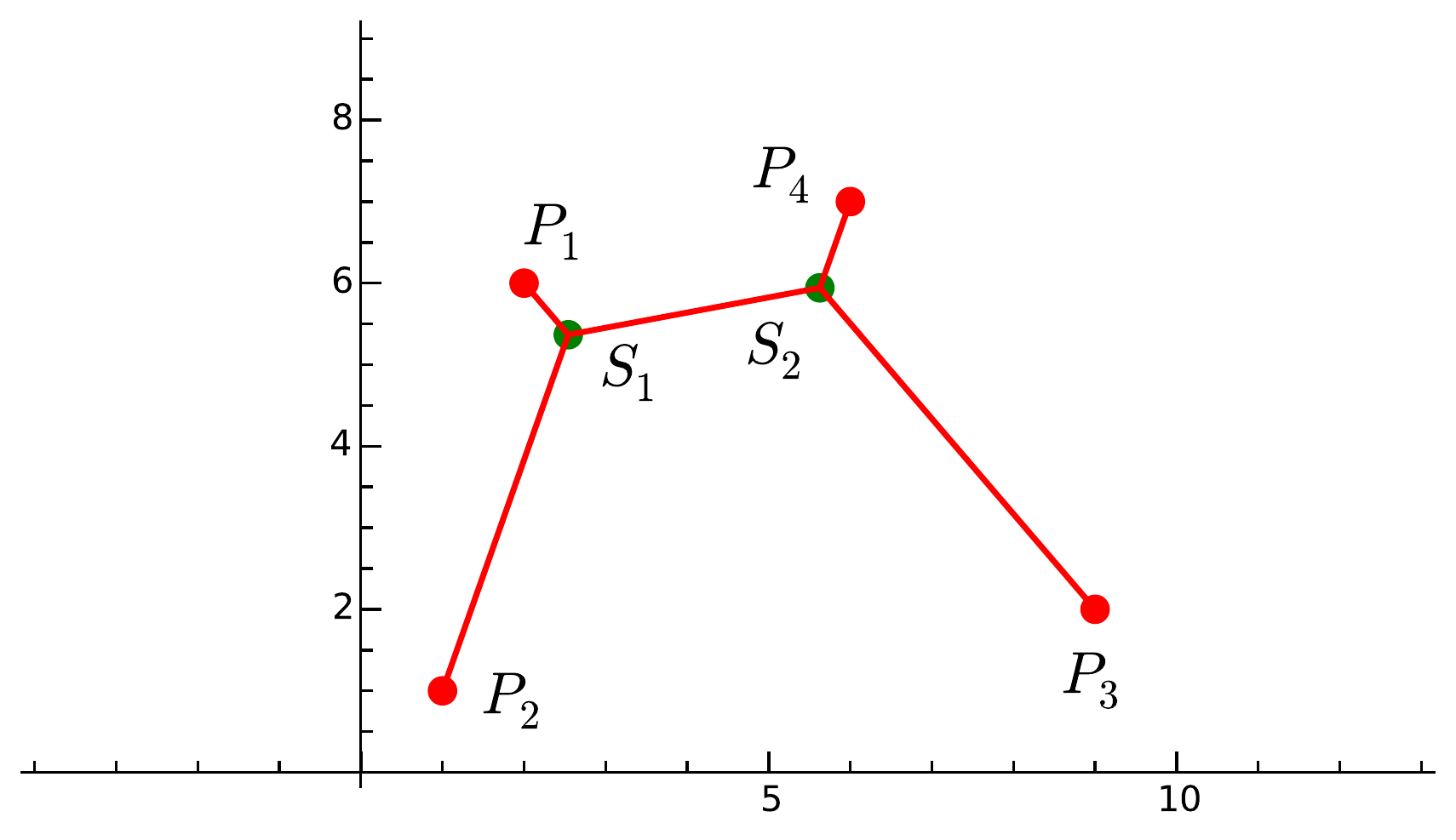}
\end{center}
\begin{center}
Figure~2. Steiner tree construction for four terminals
\end{center}
\end{figure*}
%The intersection points of the segment with the circles are Steiner points for the tree.
Coordinates of $ Q_2 $ are computed similar to $ Q_1 $:
$$
Q_2:=(q_{2x},q_{2y}) \mapsto \mathfrak q_2=q_{2x}+ \mathbf i q_{2y}=\upsilon_1z_3+\upsilon_5z_4 \, .
$$
The length of the tree equals
$$
\mathfrak L=|Q_1Q_2|=|\mathfrak q_1- \mathfrak q_2 |=|\upsilon_1z_1+\upsilon_5z_2+\upsilon_4z_3+\upsilon_2z_4|
\stackrel{(\ref{Roots_3})}{=}
| \varepsilon_2(z_3-z_1)+ \varepsilon_1(z_4-z_2)  |
$$
From this immediately follows equivalent equalities
\begin{equation}
\mathfrak L=\mathfrak L \cdot |- \varepsilon_1|=| (z_1-z_3)+ (z_2-z_4) \varepsilon_2 |=| (z_2-z_4)+ (z_1-z_3) \varepsilon_1 | \, .
\label{St41}
\end{equation}
In the case of existence a Steiner tree in the alternative topology $ \begin{array}{c} P_4 \\ P_1 \end{array} \tilde S_1 \tilde S_2 \begin{array}{c} P_3 \\ P_2 \end{array} $, its length is given by the formula
$$
\tilde \mathfrak L=| (z_1-z_3)+ (z_4-z_2) \varepsilon_1 | \, .
$$

\begin{example} \label{Ex1}  For the terminals
$$
P_1=(2,6), \ P_2= (1,1), P_3=(9,2), P_4=(6,7) \, ,
$$
and the Steiner tree topology $ \begin{array}{c} P_1 \\ P_2 \end{array} S_1 S_2 \begin{array}{c} P_4 \\ P_3 \end{array} $ (Fig.~2),
formula (\ref{St41}) yields
$$ \mathfrak L =|-7+4 \mathbf i +(-5-6 \mathbf i)\varepsilon_2 | = |-9/2-3\sqrt{3}+ \mathbf i (7+5\sqrt{3}/2)   |=\sqrt{115+62\sqrt{3}} \approx 14.912651\, .$$
\textbf{Check.} The Steiner points:
$$S_1=\left(\frac{5587}{3386}+\frac{1743}{3386}\sqrt{3}, \frac{11183}{3386}+\frac{12107}{10158}\sqrt{3}\right)\approx (2.541631, 5.367094) \, , $$
$$ S_2=\left(\frac{25479}{3386}-\frac{3711}{3386}\sqrt{3}, \frac{16193}{3386}+\frac{2267}{3386}\sqrt{3}\right)\approx  (5.626509, 5.941984) \, . $$
\end{example}

Our next aim is to discover the correlation between the sequences of terminals $ \{z_j \} $ and the roots of unity (\ref{Roots_6}) in the general formula for the  tree
length.

\section{Length of Steiner tree} \label{SLength}
\setcounter{equation}{0}
\setcounter{theorem}{0}
\setcounter{example}{0}

Consider a full Steiner tree $ \mathbb T $ with $ n \ge 3 $ terminals. Let the terminal $ P_j $ be incident to the edge $ S P_j $ linking it to the
adjacent Steiner point $ S $ of the tree. We treat the edge $ S P_j $ to be oriented from $ S $ to $ P_j $ and denote by
\begin{equation}
 \overrightarrow{\ell_j} := \overrightarrow{S P_j}/|SP_j|
 \label{direction}
\end{equation}
the \textbf{direction of the terminal} $ P_j $ \textbf{in the tree} $ \mathbb T $, i.e. a unit vector with its starting point attached to the origin.  It is known that for the full Steiner tree, the angle between $ \overrightarrow{\ell_j} $ and $ \overrightarrow{\ell_k} $ is an integer multiple of $ \pi/3 $. Therefore, all the terminals can be sorted according to the $ 6 $ possible directions. For our purpose, we do not need the exact values of these directions, but only their relative positions.

%with respect to some specified one, for instance, $ P_1 $.

\begin{theorem} \label{Th St Length} For the full Steiner tree $ \mathbb T $ with $ n\ge 3 $ terminals,  assign to each $ P_j $ the value $ U_j $ of $6$th root of unity (\ref{Roots_6})  according to the following rules:
\begin{itemize}
\item[(a)] $ P_1 \to U_1:=\upsilon_0=1 $;
\item[(b)] $ P_j \to U_j:=\upsilon_k $ if the angle between $ \overrightarrow{\ell_1} $ and $ \overrightarrow{\ell_j} $ counted \underline{clockwise} equals $ k\pi/3 $.
\end{itemize}
Then
\begin{equation}
\mathfrak L(\mathbb T)=\left| \sum_{j=1}^n z_j U_j \right| \, .
\label{Fdistance}
\end{equation}
\end{theorem}

\textbf{Proof} is by induction on $ n $. For $ n= 3 $ the statement holds as it is demonstrated in Section \ref{SHints}. Assume the
result to be true for an arbitrary full Steiner tree with $ n-1\ge 3 $ terminals. Next, consider a given full Steiner tree $ \mathbb T_n $ with $ n $ terminals, and number these terminals and Steiner points somehow with the only requirement that the terminals $ P_{n-1} $ and  $ P_{n} $ should be adjacent to the same Steiner point $ S_{n-2} $ which, in turn, should be adjacent to the Steiner point $ S_{n-3} $.

Replace the pair $ P_{n-1}, P_{n} $ by a point defined by either of alternatives
\begin{equation}
Q \mapsto \mathfrak q=\left\{ \begin{array}{c}
\upsilon_1 z_{n-1}+\upsilon_5z_n, \\
\upsilon_5 z_{n-1}+\upsilon_1z_n
\end{array}
\right.
\label{QalterN}
\end{equation}
in accordance with its location on the side of the line $ P_{n-1}P_n $ other that of $ S_{n-2} $ (Fig. 3).

\begin{minipage}[t]{80mm}
\begin{center}
\includegraphics[width=80mm]{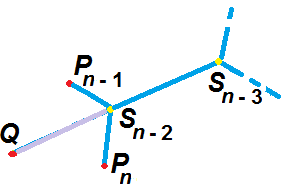}

(a)
\end{center}
\end{minipage}
%\hfill
\begin{minipage}[t]{80mm}
\begin{center}
\includegraphics[width=80mm]{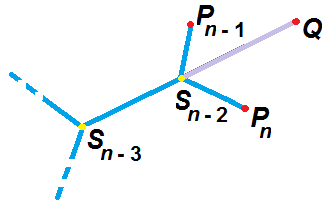}

(b)
\end{center}
\end{minipage}

\begin{center}
Figure~3. Location of the point $ Q $.
\end{center}

The new tree with terminals $ P_1,\dots,P_{n-2},Q $ and
Steiner points $ S_1,\dots, S_{n-3} $ is the full Steiner tree (the points $ S_{n-3}, S_{n-2} $ and $ Q $ are collinear). Due to Theorem \ref{ThPtolemy}, one has:
$$ |S_{n-3} S_{n-2}|+|P_{n-1}S_{n-2}|+|P_{n}S_{n-2}|=|S_{n-3}S_{n-2}| + |S_{n-2} Q |=|S_{n-3}Q| \, . $$
Therefore, the new $ (n-1) $-terminal tree $ \mathbb T_{n-1} $ has its length equal to $ \mathfrak L(\mathbb T_n) $.
The direction of the terminal $ Q $ in the tree $ \mathbb T_{n-1} $ (i.e. $ \overrightarrow{\ell}= \overrightarrow{S_{n-3} Q}/|S_{n-3} Q| $)  obeys the assignment rule from the statement of the theorem, i.e. it differs from $ \overrightarrow{\ell_1} $ by an integer multiple of $ \pi/3 $. Based on the inductive  hypothesis, one may conclude that
$$
\mathfrak L(\mathbb T_{n-1})=\left| \sum_{j=1}^{n-2} z_j U_j + \mathfrak q U \right|
$$
where the sequence $U_1,\dots,U_{n-2},U $ is computed in accordance with the assignment rule. Assume, for definiteness, that
\begin{itemize}
  \item  $ \overrightarrow{\ell} $ coincides with $ \overrightarrow{\ell_1} $, i.e. $ U=\upsilon_0=1 $,
  \item  and  $ \mathfrak q $ be defined by the first alternative from (\ref{QalterN}).
\end{itemize}
Then
$$
\mathfrak L(\mathbb T_{n})=\mathfrak L(\mathbb T_{n-1})=\left| \sum_{j=1}^{n-2} z_j U_j +  \upsilon_1 z_{n-1}+\upsilon_5z_n \right| \, .
$$
In the obtained formula, the assignment rule for terminal directions is kept fulfilled.
The direction $\overrightarrow{\ell}_{n-1} $ is $\pi/3 $ clockwise while that of $\overrightarrow{\ell_n} $ is $\pi/3 $ counterclockwise from  $\overrightarrow{\ell} $ (Fig. 3 (a)). \qed

\begin{cor} \label{Cor1} In the notation of Theorem \ref{Th St Length}, one has
$$ \sum_{j=1}^n U_j = 0 \, . $$
\end{cor}

The proof can be extracted from that of the theorem. From this statement and from the homogeneity of  (\ref{Fdistance}) with respect to the terminal coordinates, the following assertion holds

\begin{cor} \label{Cor2} The value (\ref{Fdistance}) is invariant under the affine transformation of the complex plane
$$
z \mapsto w_0z+z_0 \ ,  \forall \{w_0,z_0\} \subset \mathbb C, |w_0|=1 \, .
$$
\end{cor}

In particular, it is possible to choose such a rotation of the plane $ \mathbb R^2 $ around the origin that, in the new coordinates, the expression
$ \sum_{j=1}^n z_j U_j $ in (\ref{Fdistance}) becomes a positive real number. In this version,  formula  (\ref{Fdistance}) transforms into the interpretation by Gilbert and Pollak  of Maxwell's theorem on the equilibrium condition of the pin-jointed rigid rods holding a
prescribed system of external forces  \cite{Gilbert&Pollak, Maxwell_1864}.
However, generically, to find this rotation for an arbitrary given set $ \mathbb P $ is not a trivial task. For its evaluation, direction of at least one terminal $ P_j $  in the tree $ \mathbb T $ should be established, i.e. the coordinates of the adjacent Steiner point in (\ref{direction}).   For Example \ref{Ex1}, this ``good'' rotation is given by
$$ w_0=\cos \varphi + \mathbf i \sin \varphi \ \quad \mbox{with} \  \varphi:=\pi+ \arcsin \left( \frac{1}{3386}\sqrt{8676625-1188486\sqrt{3}} \right) \, . $$

To complete the present section, we convert (\ref{Fdistance}) into real numbers.

\begin{cor} Let $ \{U_{jx}:=\mathfrak{Re} (U_j),\ U_{jy} := \mathfrak{Im} (U_j)\}_{j=1}^n $.  Then
$$ \{U_{jx}\}_{j=1}^n \subset \{1,\pm 1/2\}, \ \{U_{jy}\}_{j=1}^n \subset \{0,\pm \sqrt{3}/2\}, \ \{U_{jx}^2+U_{jy}^2=1 \}_{j=1}^n \ . $$
Set
$$ X:=(x_1,\dots,x_n),\  Y:=(y_1,\dots,y_n), \mathfrak U_X:=(U_{1x},\dots, U_{nx}),\ \mathfrak U_Y:=(U_{1y},\dots, U_{ny}) \, .
$$
Then
$$
\mathfrak L(\mathbb T)=\sqrt{\left(\langle \mathfrak U_X,X \rangle - \langle \mathfrak U_Y,Y \rangle  \right)^2+
\left(\langle \mathfrak U_Y,X \rangle + \langle \mathfrak U_X,Y \rangle  \right)^2}
$$
where $ \langle \cdot \rangle $ is the standard inner product in $ \mathbb R^n $.
\end{cor}

\section{Examples}\label{SExamp}

We present here two examples illuminating the usage of formula (\ref{Fdistance}).

\begin{example} \label{ExSt5}  For terminals
$$ P_1=(3,9),\ P_2=(1,6),\ P_3=(6,3),\ P_4=(10,7),\ P_5=(8,10) $$
the Steiner tree is displayed in Fig. 4 (a). The assignment scheme from Theorem \ref{Th St Length} is given in Fig.~4~(b).
$$
\mathfrak L =|z_1+z_5+ \upsilon_2 z_4+\upsilon_3 z_3+\upsilon_4 z_2|=|(z_1-z_3+z_5)+\varepsilon_2z_2+\varepsilon_1z_4|
$$
$$
=\left|-\frac{1+\sqrt{3}}{2}+\mathbf i \frac{19+9\sqrt{3}}2 \right|= \sqrt{152+86\sqrt{3}} \approx 17.348094 \, .
$$
\textbf{Check.} The Steiner points
$$
\begin{array}{c|c|c}
S_1 & S_2 & S_3  \\
=\left(\frac{602}{229}+\frac{391}{1374}\sqrt{3}, \frac{2847}{458}+\frac{493}{687}\sqrt{3} \right) & = \left(\frac{1430}{229}-\frac{395}{1374}\sqrt{3}, \frac{2901}{458}-\frac{50}{687}\sqrt{3} \right) & = \left(\frac{1755}{229}+\frac{200}{687}\sqrt{3},\frac{1598}{229}+\frac{355}{687}\sqrt{3}\right)  \\
\approx ( 3.121712, 7.459099 ) & \approx (5.746608, 6.20800212  ) & \approx (8.167992,\, 7.873185)
\end{array}
$$
\end{example}

\begin{minipage}[t]{80mm}
\begin{center}
\includegraphics[width=80mm]{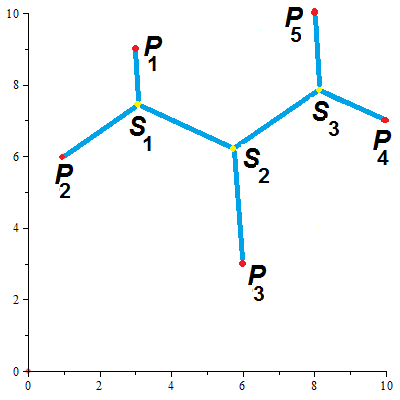}

(a)
\end{center}
\end{minipage}
%\hfill
\begin{minipage}[t]{80mm}
\begin{center}
\includegraphics[width=50mm]{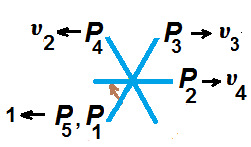}

   (b)
\end{center}
\end{minipage}

\begin{center}
Figure~4. (a) The Steiner tree for $ 5 $ terminals and (b) its assignment scheme (Example \ref{ExSt5}).
\end{center}

\begin{example} \label{ExSt6} For terminals
$$P_1=(2,-2\sqrt{3}), P_2=(8,0), P_3=(4,4\sqrt{3}), P_4=(-8/\sqrt{3}, 8), P_5= (-10,0), P_6=(-2 , -2\sqrt{3})
$$
the Steiner tree is displayed in Fig.~5~(a). The assignment scheme is given in Fig.~5~(b)
$$
\mathfrak L=|z_1+\upsilon_1 z_5+ \upsilon_2 z_4+\upsilon_3 z_3+\upsilon_4 z_2+\upsilon_5 z_6|=|(z_1-z_3)+\varepsilon_1(z_4-z_6)+\varepsilon_2(z_2-z_5)|
$$
$$
=|(-15-8/\sqrt{3})- \mathbf i (8+15\sqrt{3}) |=30+15/\sqrt{3} \approx 39.237604 .
$$
\textbf{Check.} The Steiner points:
$$  S_1=(4,0), \ S_2=(2,2\sqrt{3}), S_3=(2,-2\sqrt{3}),\ S_4= (-4,0) \, . $$
\end{example}

\begin{minipage}[t]{80mm}
\begin{center}
\includegraphics[width=80mm]{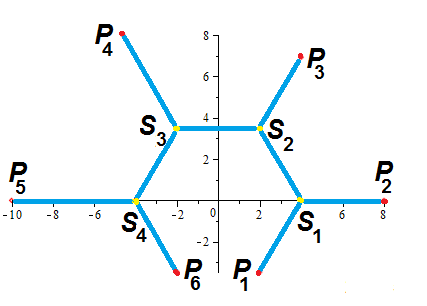}

(a)
\end{center}
\end{minipage}
%\hfill
\begin{minipage}[t]{80mm}
\begin{center}
\includegraphics[width=40mm]{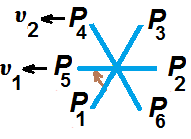}

   (b)
\end{center}
\end{minipage}

\begin{center}
Figure~5. (a) The Steiner tree for $ 6 $ terminals and (b) its assignment scheme (Example \ref{ExSt6}).
\end{center}

\section{Extension to the Weber problem} \label{SWeber_L}
\setcounter{equation}{0}
\setcounter{theorem}{0}
\setcounter{example}{0}

We now turn to the Weber problem mentioned in Introduction. In the present section we assume the fulfillment of the conditions for the existence of the networks in the prescribed
topologies (some of them are presented in \cite{Uteshev_Semenova_2020_Ar,Uteshev_Semenova_2020_LNCS}). In particular, we suppose all the involved geometric constructions to be feasible.
Hereinafter all the weights $ m_j, m, \dots $ are positive real numbers and the angles denoted by Greek letters $ \alpha, \beta, \gamma $ are assumed to be within the interval $ [0, \pi] $.

We start with the three terminal problem (also known as the generalized Fermat-Torricelli problem) as that of finding the facility $ W $ providing
\begin{equation}
\min_{W\in \mathbb R^2} \sum_{j=1}^3 m_j |WP_j| \, .
\label{Weber_uni}
\end{equation}
A geometric solution to the problem is given in the paper by Launhardt \cite{Dingeldey,Launhardt}; it is just a counterpart of the Torricelli-Simpson construction for the Steiner problem from Section \ref{SHints}.

\begin{example}\label{Exft0}
Find the optimal position of the facility $ W $ to the problem (\ref{Weber_uni}) where
$$
\left\{\begin{array}{c|c|c}
P_1=(2,6) & P_2=(1,1) & P_3=(5,1)  \\
m_1=2 & m_2=3 & m_3=4
\end{array} \right\}.
$$
\end{example}

\textbf{Solution.} First find the point $ Q_1 $ lying on the opposite side of the line $ P_1P_2 $ with respect to the point $ P_3 $ and such that
\begin{equation}
 |P_1Q_1|=\frac{m_2}{m_3}|P_1P_2|, \ |P_2Q_1|=\frac{m_1}{m_3}|P_1P_2| \, .
 \label{SimTr}
\end{equation}
These relations mean that the triangle $ P_1P_2Q_1 $ is similar to the so-called \textbf{weight triangle} of the problem, i.e. the triangle containing the edges formally coinciding with the values of the weights $ m_1, m_2, m_3 $. We will further denote the angles of the weight triangle by $ \alpha_1,\alpha_2, \alpha_3 $ as shown in Fig.~6 (a).
Next, draw the circle $ C_1 $ circumscribing $P_1P_2Q_1$. Finally draw the line through $ Q_1 $ and $ P_3 $. The intersection point of this line with $ C_1 $ is the position of the optimal facility $ W $.
It possesses the property that directions of the terminals, i.e.  $ \overrightarrow{\ell_j} := \overrightarrow{W P_j}/|WP_j| $, make the following angles (Fig.~6~(b))
$$ \widehat{\overrightarrow{\ell_1}, \overrightarrow{\ell_2}} = \pi - \alpha_3, \ \widehat{\overrightarrow{\ell_1}, \overrightarrow{\ell_3}} = \pi - \alpha_2, \ \widehat{\overrightarrow{\ell_2}, \overrightarrow{\ell_3}} = \pi - \alpha_1 \, .
$$
The corresponding (minimal) cost equals $ m_3|P_3Q_1| $ (consequence of (\ref{SimTr}) and Theorem \ref{ThPtolemy}).  \qed

\begin{minipage}[t]{80mm}
\begin{center}
\includegraphics[width=80mm]{{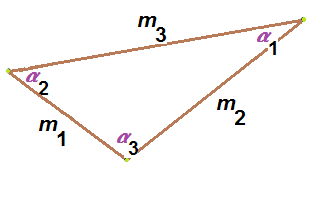}}

(a)
\end{center}
\end{minipage}
%\hfill
\begin{minipage}[t]{80mm}
\begin{center}
\includegraphics[width=80mm]{{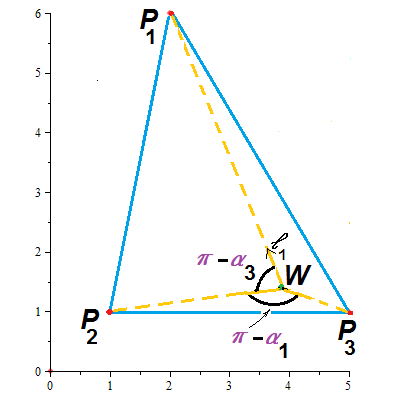}}

   (b)
\end{center}
\end{minipage}

\begin{center}
Figure~6. (a) Weight triangle and (b) Weber point location (Example \ref{Exft0}).
\end{center}

We now rewrite the analytics similar to that from Section \ref{SHints}. The coordinates of the point $ Q_1 $ are as follows \cite{Uteshev_Semenova_2020_Ar}
\begin{eqnarray}
Q_1&=&\Bigg(\frac{1}{2}(x_1+x_2)+\frac{(m_1^2-m_2^2)(x_1-x_2)-\sqrt{\mathbf k}(y_1-y_2)}{2m_3^2}, \nonumber \\
    & &\hspace{20mm}\frac{1}{2}(y_1+y_2)+\frac{(m_1^2-m_2^2)(y_1-y_2)+\sqrt{\mathbf k}(x_1-x_2)}{2m_3^2} \Bigg) \, ,
\label{Q1}
\end{eqnarray}
or, equivalently, in the complex plane:
\begin{equation}
Q_1 \mapsto \mathfrak q_1= \frac{1}{2m_3^2} \left[ (m_3^2+m_1^2-m_2^2 + \mathbf i \sqrt{\mathbf k}) z_1 + (m_3^2-m_1^2+m_2^2 - \mathbf i \sqrt{\mathbf k}) z_2  \right]
\label{Q1_alt}
\end{equation}
Here
\begin{equation}
\mathbf k:= (m_1+m_2+m_3)(-m_1+m_2+m_3)(m_1-m_2+m_3)(m_1+m_2-m_3) \, ,
\label{k}
\end{equation}
and, due to Heron's formula, $ \sqrt{\mathbf k}/4 $ equals the square of the weight triangle. As for the angles of this triangle, by cosine theorem, one gets
\begin{equation}
\frac{m_3^2+m_1^2-m_2^2}{2m_3m_1} = \cos \alpha_2, \  \frac{m_3^2-m_1^2+m_2^2}{2m_3m_2} = \cos \alpha_1  ,
\label{cosine}
\end{equation}
while\footnote{We recall the assumption made at the beginning of the present section: all the angles denoted by Greek letters are assumed to lie within $ [0, \pi] $. Therefore all the sine values for such angles are nonnegative.}
\begin{equation}
\frac{\sqrt{\mathbf k}}{2 m_3 m_1} = \sin \alpha_2,\ \frac{\sqrt{\mathbf k}}{2 m_3 m_2} = \sin \alpha_1 \, .
\label{sine}
\end{equation}
Finally, one has
$$
\mathfrak q_1= \frac{1}{m_3} \left[m_1 (\cos \alpha_2 + \mathbf i \sin \alpha_2)z_1+m_2 (\cos \alpha_1 - \mathbf i \sin \alpha_1)z_2 \right]
$$
and
$$
\mathfrak C =m_3 |P_3Q_1| =|m_3z_3 - m_3 \mathfrak q_1|
$$
$$
=| m_3z_3 +m_2z_2 \left\{\cos (\pi-\alpha_1) + \mathbf i \sin (\pi-\alpha_1)\right\} +m_1z_1 \left\{\cos (2\pi-\alpha_1-\alpha_3) + \mathbf i \sin (2\pi-\alpha_1-\alpha_3)\right\} | \, .
$$
The obtained formula is the true counterpart of the formula (\ref{Length3}) for the length of the full Steiner tree  for $3 $ terminals.

Let us now turn to the bifacility Weber problem for $ 4 $ terminals $ \{P_j\}_{j=1}^4 $
as that of finding the points $ W_1 $ and $ W_2 $ which yield
\begin{eqnarray}
& \displaystyle \min_{\{W_1,W_2\} \subset \mathbb R^2} F(W_1,W_2)  \quad \mbox{ where } & \nonumber \\
& F(W_1,W_2)= m_1|W_1P_1|+m_2|W_1P_2| +m_3|W_2P_3|+m_4|W_2P_4|+ m |W_1W_2| \, .
\label{F_Weber}
\end{eqnarray}
We  first recall the geometric solution outlined
by Georg Pick in the Mathematical Appendix of Weber's book \cite{Weber}.
We illustrate this algorithm with the following example.

\begin{example}\label{ExWeber4}
Find the optimal position for the facilities $ W_1 $ and $ W_2 $ for the problem (\ref{F_Weber}) where
$$
\left\{\begin{array}{c|c|c|c|}
P_1=(1,5) & P_2=(2,1) & P_3=(7,2) & P_4=(6,7) \\
m_1=3 & m_2=2 & m_3=3 & m_4=4
\end{array} %\mbox{ and }
 \ m=4  \right\} \, .
$$
\end{example}

\textbf{Solution.} First find the point $ Q_1 $ lying on the opposite side of the line $ P_1P_2 $ with respect to the point $ P_3 $  and such that
\begin{equation}
 |P_1Q_1|=\frac{m_2}{m}|P_1P_2|, \ |P_2Q_1|=\frac{m_1}{m}|P_1P_2| \, .
 \label{similarity1}
\end{equation}
The exact coordinates of this point are given by (\ref{Q1}) where the substitution $ m_3 \to m $ is made.
Find then the second point $ Q_2 $ with the similar property with respect to the points $ P_3 $ and $ P_4 $ (Fig. 7):
$$ |P_3Q_2|=\frac{m_4}{m}|P_3P_4|, \ |P_4Q_2|=\frac{m_3}{m}|P_3P_4| \, . $$

\begin{center}
\begin{minipage}[t]{125mm}
\includegraphics[scale=0.3]{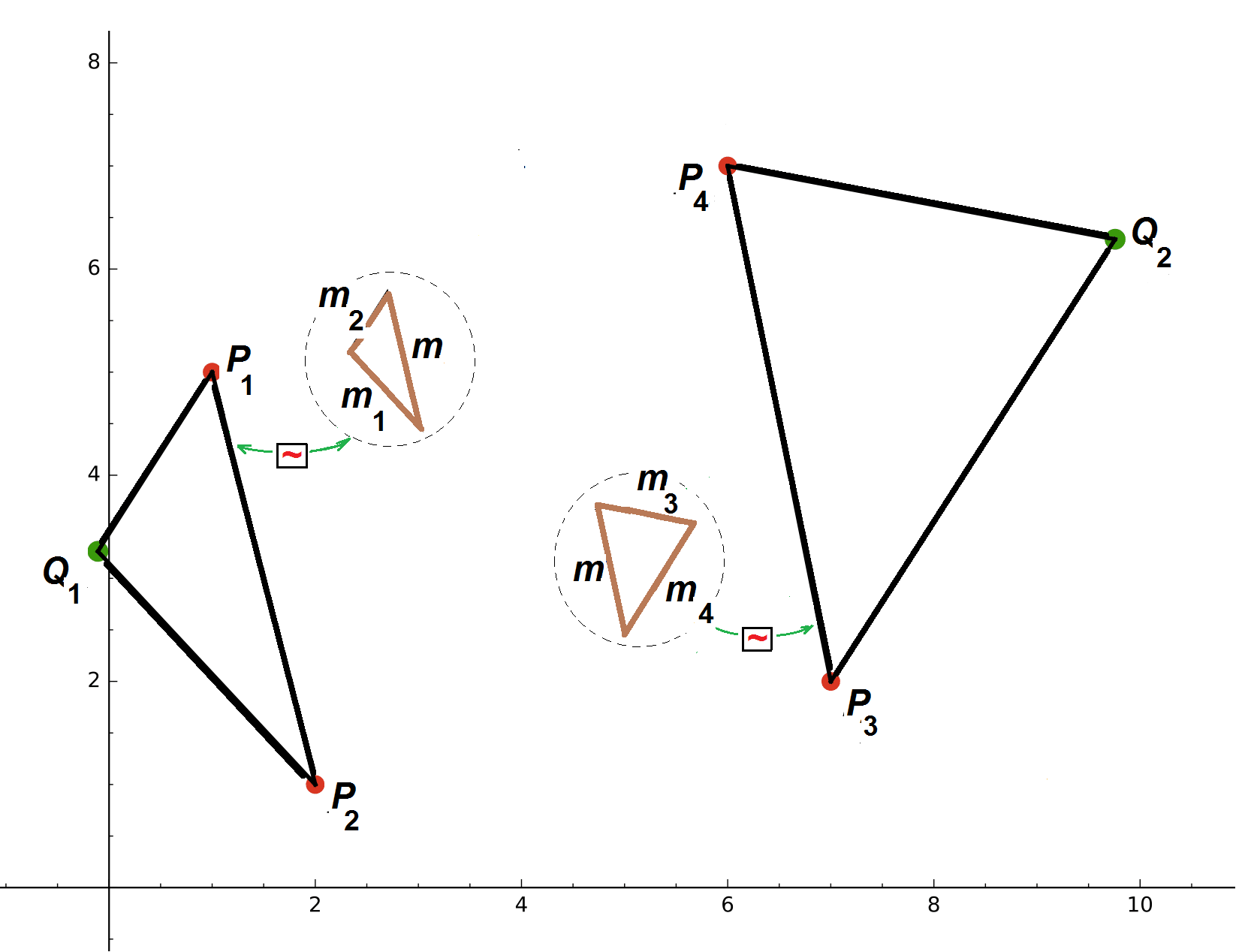}
\end{minipage}
\end{center}
%\caption{}
\begin{center}
Figure~7. Construction of the points $ Q_1 $ and $ Q_2 $.
\end{center}
It is evident that triangle $ P_1P_2Q_1 $ is similar to the weight triangle with the edges $ m, m_1,m_2 $, while $P_3P_4Q_2 $ is similar to the weight triangle with the edges $ m, m_3, m_4 $ (Fig. 8)
Next, draw the circle $ C_1 $ circumscribing $P_1P_2Q_1$ and $ C_2 $ circumscribing $P_3P_4Q_2$.

\begin{center}
\begin{minipage}[t]{125mm}
\includegraphics[scale=0.5]{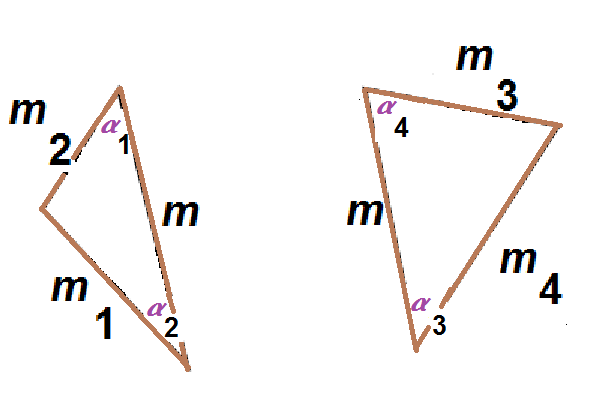}
\end{minipage}
\end{center}
%\caption{}
\begin{center}
Figure~8. Weight triangles for the bifacility Weber problem
\end{center}

Finally draw the line through $ Q_1 $ and $ Q_2 $ (Fig. 9).

{
\begin{minipage}[t]{145mm}
\begin{minipage}[t]{85mm}
%\graphicspath{{Illustrations/}}
\includegraphics[scale=0.4]{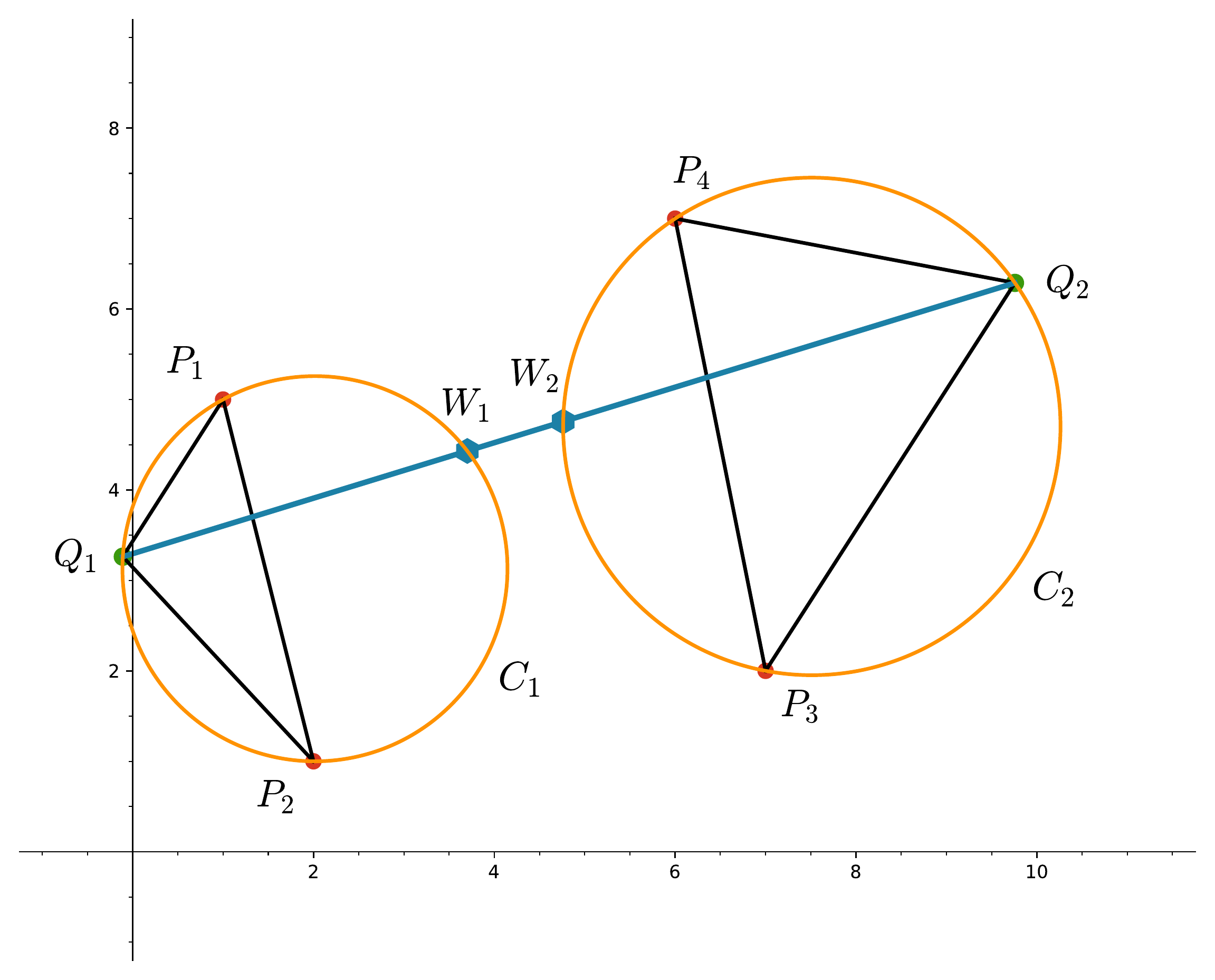}
\end{minipage}
\hfill
\begin{minipage}[t]{60mm}
%\graphicspath{{Illustrations/}}
%\includegraphics[width=60mm]{WeightQuadr.pdf}
\includegraphics[scale=0.4]{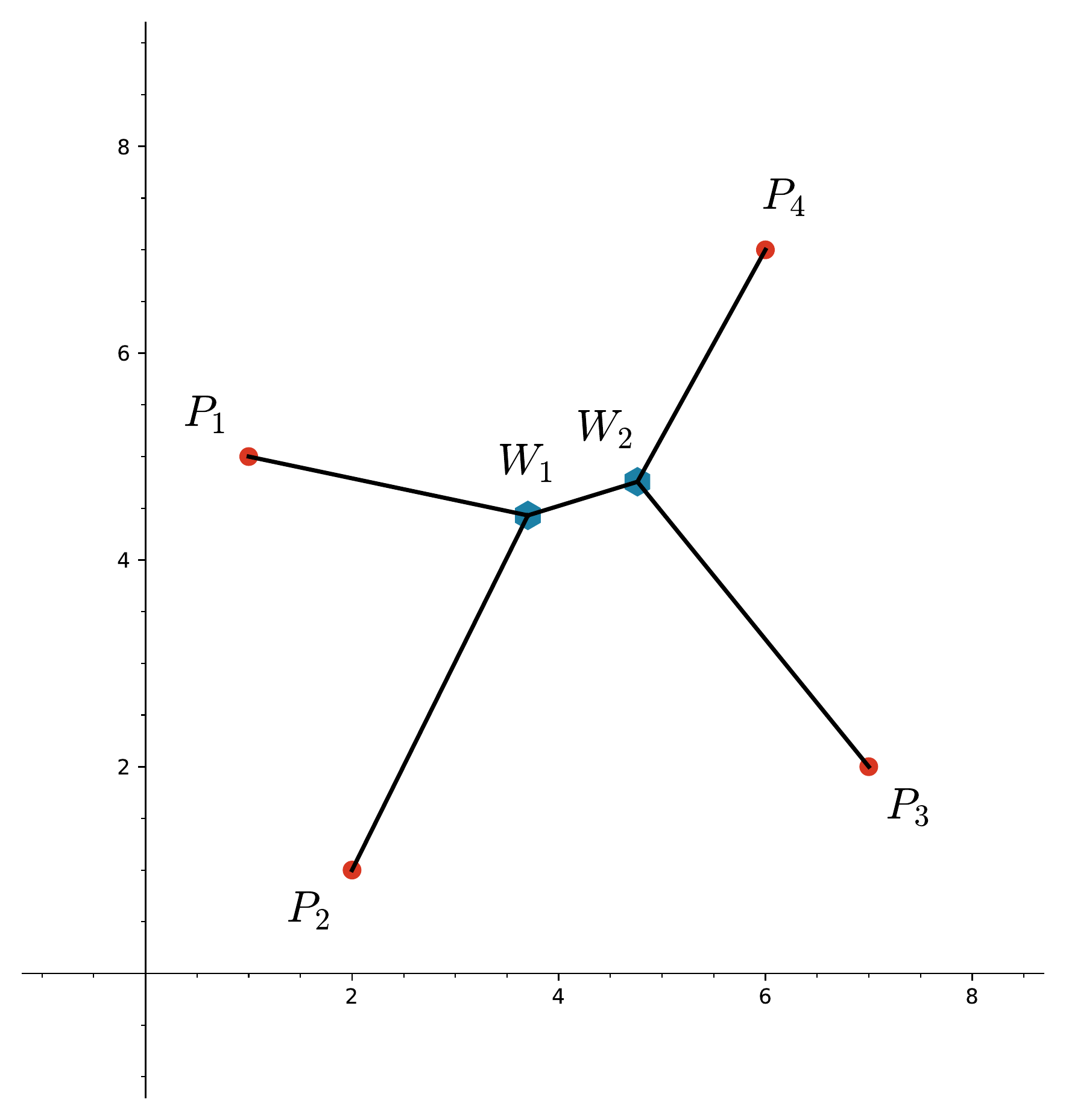}
\end{minipage}
\end{minipage}
}
%\caption{Example \ref{ft}. Pick's construction of the Weber network}  \label{Pick_4}

\begin{center}
Figure~9.  Pick's construction for the Weber network (Example \ref{ExWeber4}).
\end{center}

The intersection points of this line with $ C_1 $ and $ C_2 $ are the position of the optimal facilities $ W_1 $ and $ W_2 $ for the network
with the corresponding (minimal) cost equal to $ \mathfrak C= m|Q_1Q_2| $.

\textbf{Remark.} Pick's solution can be interpreted as a counterpart of the algorithm worked out  by Gergonne in 1810  (and rediscovered by Melzak in 1961) for solution of the SMT problem for four terminals. It should be noted that Pick did not provide any reasons of validity for his construction (we also failed to find any references to Pick's solution in subsequent papers on the subject). In \cite{Uteshev_Semenova_2020_Ar, Uteshev_Semenova_2020_LNCS} we have verified the accuracy of Pick's treatment via an analytical
representation of the coordinates of the optimal facilities $ W_1, W_2 $ and substitution them into the gradient of the function (\ref{F_Weber}).
Some other qualitative conclusions have been deduced, such as, for instance, that, in the case of existence, the bifacility network is less costly than
any unifacility one. All the computations for the examples below can be found in the cited sources.

Thus, the coordinates for the facilities from the previous example  are as follows
\begin{eqnarray*}
W_1 &=& \Bigg(\frac{\scriptstyle{2266800+772027 \sqrt{15}+453552 \sqrt{33}+246177 \sqrt{55}}}{\scriptstyle{48 \left(22049+2085 \sqrt{15}+945
   \sqrt{33}+2559 \sqrt{55}\right)}},
    \frac{\scriptstyle{1379951+201984 \sqrt{15}+97279 \sqrt{33}+154368 \sqrt{55}}}{\scriptstyle{16 \left(22049+2085 \sqrt{15}+945
   \sqrt{33}+2559 \sqrt{55}\right)}} \Bigg) \\
   & & \approx ( 3.701271,  4.430843) \, ; \\
W_2 &=& \Bigg(\frac{\scriptstyle{188467345+18613485\sqrt{15}+7149825\sqrt{33}+20949207\sqrt{55}}}{\scriptstyle{1760\left(22049+2085 \sqrt{15}+945
   \sqrt{33}+2559 \sqrt{55}\right)}},  \frac{\scriptstyle{188346565+19265895\sqrt{15}+20525157\sqrt{55}+7187445 \ \sqrt{11}}}
   {\scriptstyle{1760\left(22049+2085 \sqrt{15}+945
   \sqrt{33}+2559 \sqrt{55}\right)}} \Bigg)  \\
   & & \approx ( 4.761622, 4.756175)
\end{eqnarray*}
and the cost of the network
\begin{equation}
\mathfrak C=\frac{1}{8}\sqrt{44098+4170\sqrt{15}+5118\sqrt{55}+1890\sqrt{33}}  \approx  41.280608 \, .
\label{CostEx}
\end{equation}

Let us deduce the general formula for $ \mathfrak C= m|Q_1Q_2| $ for the networks which are topologically equivalent to the
one dealt with in Example \ref{ExWeber4}.
The coordinates of the point $ Q_1 $ are determined by (\ref{Q1_alt}) where substitution $ m_3 \to m $ is made. The coordinates of $ Q_2 $ are also obtained from (\ref{Q1_alt})
by replacement $ (m_1,m_2,m_3) \to (m_3,m_4, m) $. Thus
\begin{eqnarray}
\mathfrak C= m|Q_1Q_2| = m |\mathfrak q_1 - \mathfrak q_2| &=&\big|m_1z_1 (\cos \alpha_2 + \mathbf i \sin \alpha_2)+m_2z_2 (\cos \alpha_1 - \mathbf i \sin \alpha_1)  \nonumber \\
& & - m_3z_3 (\cos \alpha_4 + \mathbf i \sin \alpha_4)-m_4z_4 (\cos \alpha_3 - \mathbf i \sin \alpha_3) \big| \label{Cost4} \, .
\end{eqnarray}
Here $ \alpha_1, \alpha_2 $ are the angles of the first weight triangle of the problem, while $ \alpha_3, \alpha_4 $ are those of the second (Fig. 8). Cosine and sine functions of these angles are computed similar to (\ref{cosine}) and (\ref{sine}).
For the network of Example \ref{ExWeber4} one gets the cost value
$$
\mathfrak C = \left| 3(1+\mathbf i 5)\left(\frac{7}8+\frac{1}{24}\mathbf i \sqrt{135} \right)+2(2+\mathbf i) \left(\frac{11}{16}-\frac{1}{16}\mathbf i \sqrt{135}\right) - \dots \right|
$$
$$
= \left| -\frac{79}4-\frac{1}{2}\sqrt{135}-\frac{5}{8}\sqrt{495} - \mathbf i \left( \frac{63}8+\frac{1}8\sqrt{135}+\frac{1}8\sqrt{495} \right) \right|
$$
which coincides with (\ref{CostEx}).

\begin{minipage}[t]{80mm}
\begin{center}
\includegraphics[width=80mm]{{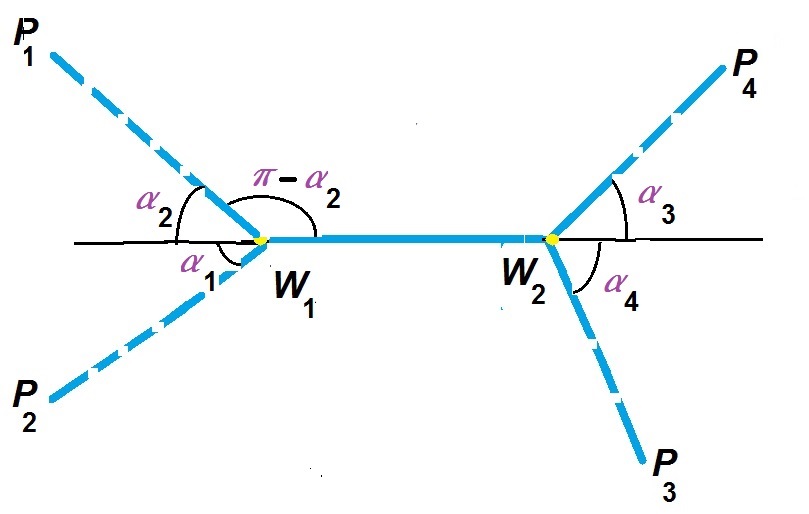}}

(a)
\end{center}
\end{minipage}
%\hfill
\begin{minipage}[t]{80mm}
\begin{center}
\includegraphics[width=50mm]{{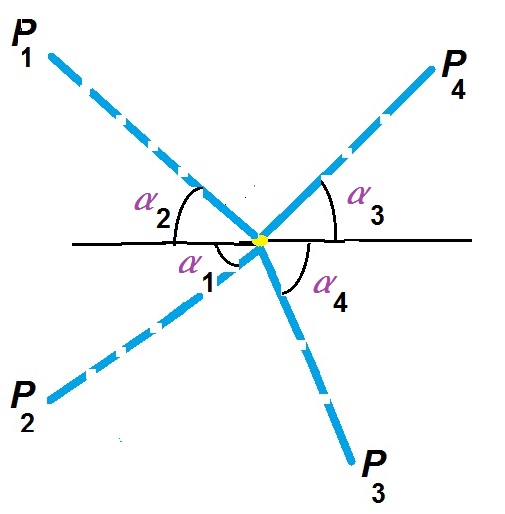}}

   (b)
\end{center}
\end{minipage}

\begin{center}
Figure~10. Angles in the bifacility Weber network.
\end{center}
Let us rewrite (\ref{Cost4}) uniformizing representation of the arguments of trigonometric functions:
\begin{eqnarray}
\mathfrak C&=&\big|m_1 (\cos \alpha_2 + \mathbf i \sin \alpha_2)z_1+m_4z_4 (\cos (\pi-\alpha_3) + \mathbf i \sin (\pi-\alpha_3)) \nonumber  \\
&& +m_3z_3 (\cos (\pi+\alpha_4) + \mathbf i \sin (\pi+ \alpha_4))
+m_2z_2 (\cos (2\pi-\alpha_1) + \mathbf i \sin (2\pi-\alpha_1)) \big| \, . \label{Cost42}
\end{eqnarray}
Here one can watch the generation rule for these arguments: they correspond to the direction angles of the terminals $ P_j $
in the network if we count them clockwise starting from the line $ W_1W_2 $ (Fig. 10).

\begin{theorem} \label{ThWeberC} Let the optimal Weber network $ \mathbb T $ with $ n \ge 3 $ terminals $ \{ P_j \}_{j=1}^n $ and $ n-2 $ facilities $ \{W_k\}_{k=1}^{n-2} $ be topologically equivalent to a full Steiner tree, i.e. each terminal from $ \mathbb T $ is of the degree $ 1 $  while each facility  is of the degree $ 3 $. Let $  \overrightarrow{\ell_j}  =(\cos \beta_j, \sin \beta_j ) $
be the direction of the terminal $ P_j $ in the network counted clockwise starting from some particular vector.
The cost of the network $ \mathbb T $ in the prescribed topology is given by the formula
\begin{equation}
\mathfrak C (\mathbb T)=\left| \sum_{j=1}^n  m_j z_j (\cos \beta_j+ \mathbf i \sin \beta_j) \right| \, .
\label{Wcost}
\end{equation}
\end{theorem}

The results of Corollaries \ref{Cor1} and \ref{Cor2} from Section \ref{SLength} can evidently be upgraded to the ``weighted'' version. For instance,
one has
\begin{equation}
\sum_{j=1}^n m_j (\cos \beta_j+ \mathbf i \sin \beta_j) = 0 \, .
\label{WcostSum}
\end{equation}

To justify the result of the theorem, we suggest here some \emph{plausible reasonings}.

\begin{example}\label{Ex5t3f} Find the cost of the optimal network that minimize the function
\begin{equation}
m_1|P_1W_1|+m_2|P_2W_1|+m_3|P_3W_2| +  m_4|P_4W_2|+m_5|P_5W_3|
    +  \widetilde{m}_{1,3} |W_1W_3| + \widetilde{m}_{2,3} |W_2W_3|
\label{cost5-3}
\end{equation}
for the following configuration:
$$
\left\{\begin{array}{c|c|c|c|c|c}
P_1=(1,6) & P_2=(5,1) & P_3=(11,1) & P_4=(15,3) & P_5=(7,11) & \widetilde{m}_{1,3}=10  \\
m_1=10 & m_2=9 & m_3=8 & m_4=7 & m_5=13 & \widetilde{m}_{2,3}=12
\end{array}
\right\} \, .
$$
\end{example}

\textbf{Solution.} Existence of the optimal network and coordinates of the corresponding facilities $ W_1, \, W_2, \, W_3 $
are established in \cite{Uteshev_Semenova_2020_Ar}. We will discuss here only the minimal cost calculation.
\begin{center}
\includegraphics[width=0.70\textwidth]{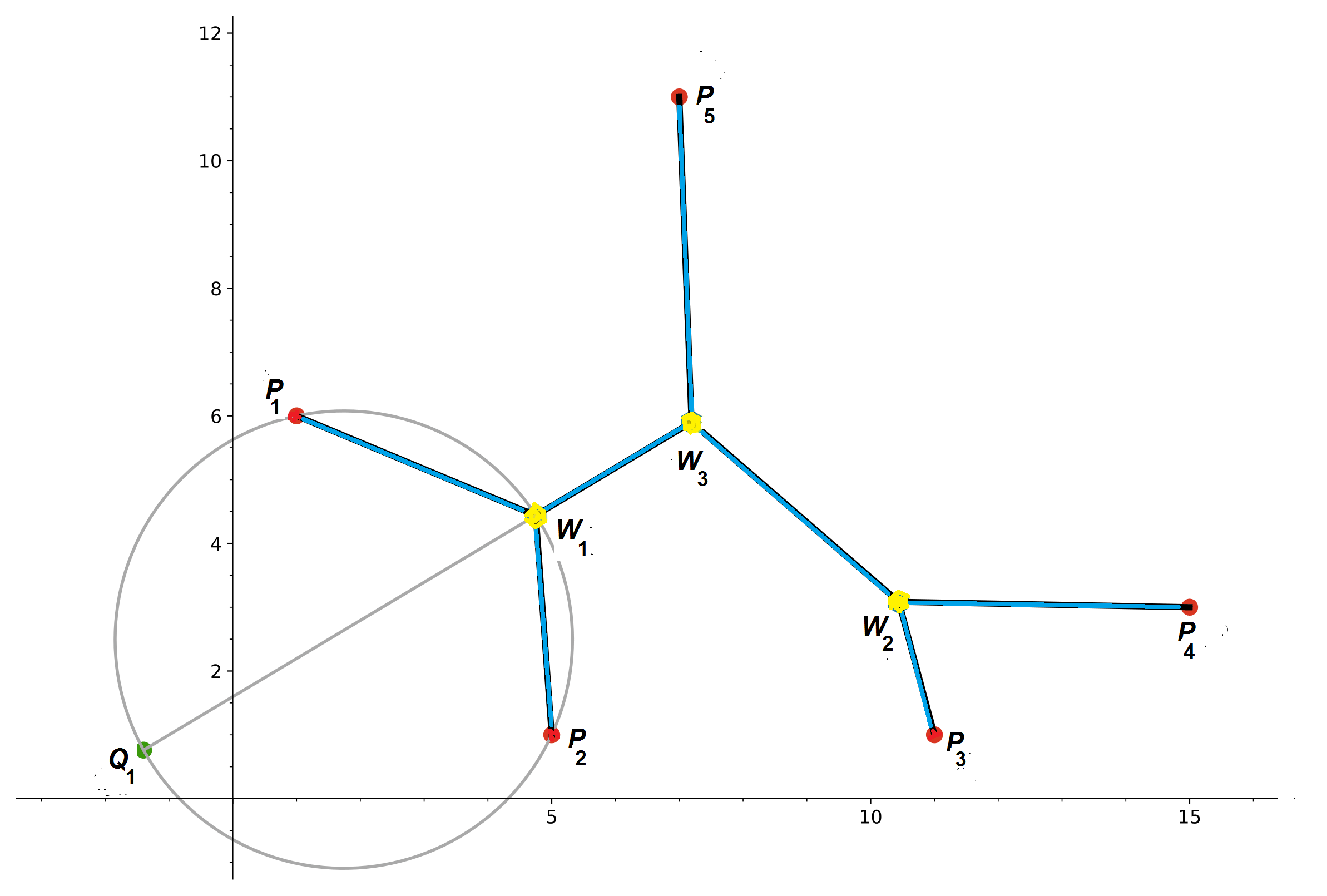}
\end{center}
\begin{center}
Figure~11. Weber network construction for five terminals (Example \ref{Ex5t3f}).
\end{center}
Replace a pair of the terminals $ P_1 $ and $ P_2 $ by the \emph{phantom terminal}
$ Q_1 $ defined by the formula (\ref{Q1_alt}) where the substitution $ m_3 \to \widetilde{m}_{1,3} $ is made, and assign this weight to
$ Q_1 $.
%$$
%Q_1=\left( -\frac{9}{40}\sqrt{319} + \frac{131}{50} , -\frac{9}{50}\sqrt{319} + \frac{159}{40} \right) \approx ( -1.398628 , 0.760097) \, .
%$$
The minimal cost of the new $\left\{ 4 \right. $-terminals, $ 2 $-facilities$\left. \right\}$-network minimizing the function
$$
\widetilde{m}_{1,3} |Q_1W_3|+ m_3|P_3W_2|+m_4|P_4W_2|+m_5 |P_5W_3|
$$
equals the minimal cost of the original network (and facilities
$ W_2, W_3 $ in both optimal networks coincide).  By (\ref{Cost42}),
\begin{eqnarray}
\mathfrak C&=&\big|\widetilde{m}_{1,3} \mathfrak q_1 (\cos (2\pi-\gamma_1) + \mathbf i \sin (2\pi-\gamma_1))+m_5z_5 (\cos \gamma_2 + \mathbf i \sin \gamma_2))  \nonumber \\
& & +m_4z_4 (\cos (\pi-\gamma_3) + \mathbf i \sin (\pi-\gamma_3))+m_3z_3 (\cos (\pi+\gamma_4) + \mathbf i \sin (\pi+\gamma_4)) \big|
\label{Cost4-9}
\end{eqnarray}
with the angles $ \{\gamma_j \}_{j=1}^4 $ displayed in Fig.~12~(a).

\begin{minipage}[t]{80mm}
\begin{center}
\includegraphics[width=80mm]{{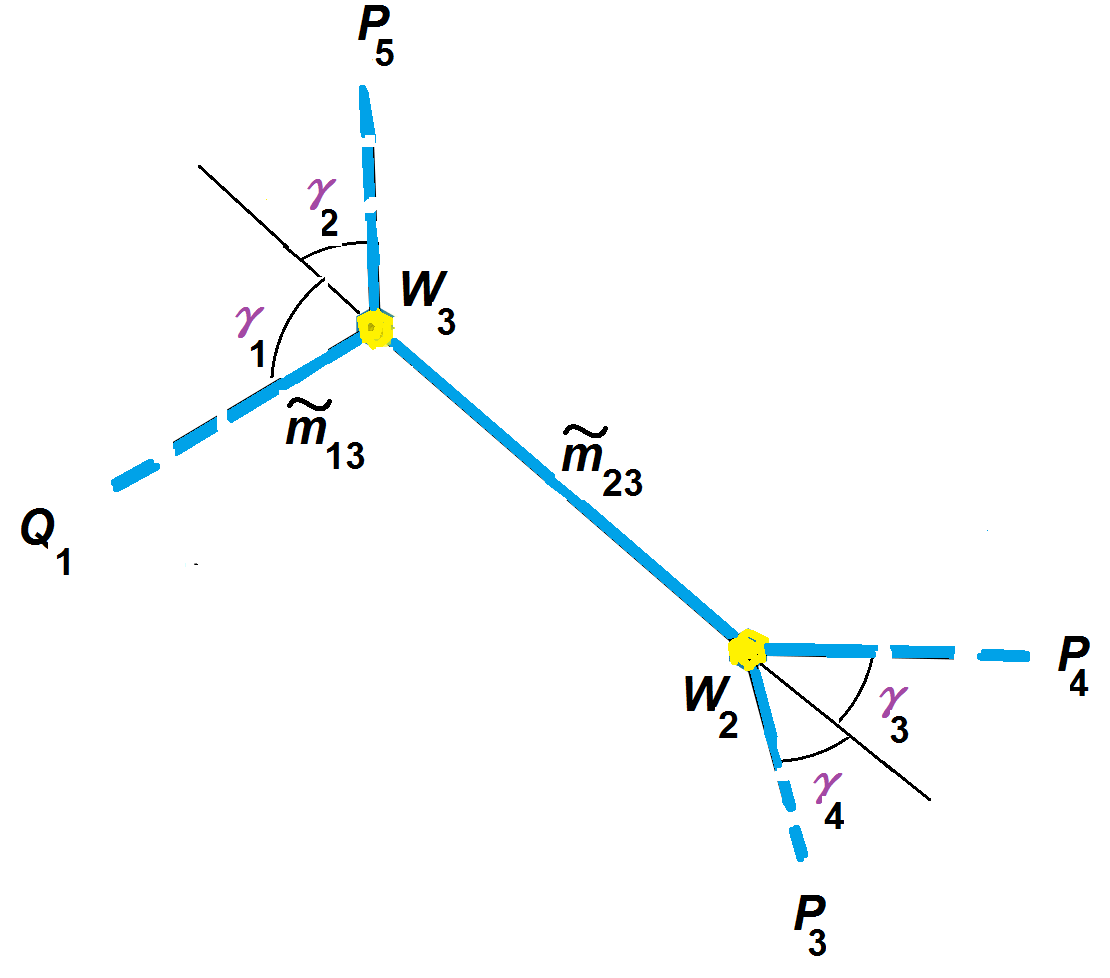}}

(a)
\end{center}
\end{minipage}
%\hfill
\begin{minipage}[t]{80mm}
\begin{center}
\includegraphics[width=80mm]{{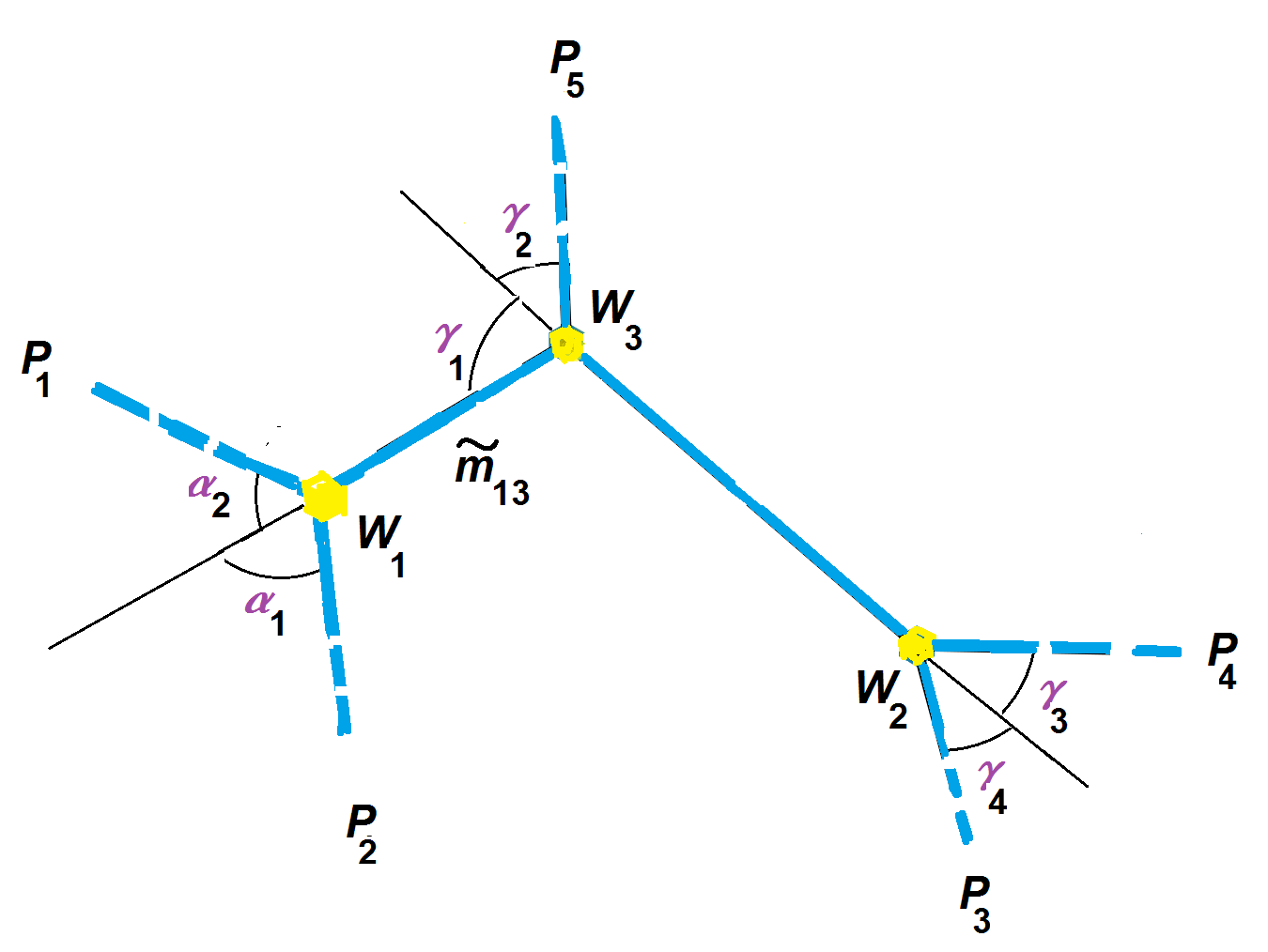}}

   (b)
\end{center}
\end{minipage}
\begin{center}
Figure~12. Angles of the trifacility Weber network (Example \ref{Ex5t3f}).
\end{center}

One has
$$
\cos \gamma_1=\frac{\widetilde{m}_{1,3}^2+\widetilde{m}_{2,3}^2-m_5^2}{2\widetilde{m}_{1,3}\widetilde{m}_{2,3}},\
\cos \gamma_2=\frac{\widetilde{m}_{2,3}^2+m_5^2-\widetilde{m}_{1,3}^2}{2m_5\widetilde{m}_{2,3}},
$$
$$
\cos \gamma_3=\frac{m_4^2+\widetilde{m}_{2,3}^2-m_4^2}{2\widetilde{m}_{2,3}m_4},\
\cos \gamma_4=\frac{m_3^2+\widetilde{m}_{2,3}^2-m_3^2}{2\widetilde{m}_{2,3}m_3}
$$
and the expressions for the corresponding sine values are represented, using the involved weights,  similarly to (\ref{sine}). Next we replace $ \mathfrak q_1$ by the original
terminals $ z_1 $ and $ z_2 $ via (\ref{Q1_alt})  ($ m_3 \to \widetilde{m}_{1,3} $):
$$
\mathfrak q_1 = \frac{1}{\widetilde{m}_{1,3}} \left[m_1 (\cos \alpha_2 + \mathbf i \sin \alpha_2)+
m_2 (\cos (2\pi-\alpha_1) + \mathbf i \sin (2\pi-\alpha_1)) \right]
$$
with the angles $ \alpha_1, \alpha_2 $ displayed in Fig. 12 (b) and with
$$
\cos \alpha_1=\frac{\widetilde{m}_{1,3}^2+m_2^2-m_1^2}{2m_2\widetilde{m}_{1,3}},\
\cos \alpha_2=\frac{\widetilde{m}_{1,3}^2+m_1^2-m_2^2}{2m_1\widetilde{m}_{1,3}}\, .
$$
Finally,
$$
\mathfrak C=\big|m_1z_1 (\cos (\alpha_2-\gamma_1) + \mathbf i \sin (\alpha_2-\gamma_1))+m_2z_2
(\cos (2\pi-\alpha_1-\gamma_1) + \mathbf i \sin (2\pi-\alpha_1-\gamma_1))
+ \dots  \big|
$$
with the rest of the terms coinciding with those from (\ref{Cost4-9}). The terminal direction assignment rule from Theorem \ref{ThWeberC} is fulfilled: these angles are  counted starting from the line $ W_2W_3 $. Substitution of the expressions for sine and cosine functions in terms
of all the involved weights yields the cost representation by radicals. For our particular example,
\begin{eqnarray*}
\mathfrak C&=& \Big|
-\frac{1331}{16}-\frac{281}{64}\sqrt{231}-\frac{45}{64}\sqrt{319}-\frac{3}4\sqrt{143} -\frac{9}{80}\sqrt{319}\sqrt{231} \\
&& + \mathbf i \left( \frac{5587}{64}+\frac{219}{80}\sqrt{231}-\frac{9}{16}\sqrt{319}+\frac{3}{2}\sqrt{143}+\frac{9}{64}\sqrt{319}\sqrt{231}\right)
\Big| \approx 267.229644 \, .
\end{eqnarray*}
\qed

The result of the theorem remains valid even for the networks  without the imposed restriction on their topologies. For instance, the unifacility Weber problem (generalized Fermat-Torricelli problem) for $ n\ge 4 $ terminals as that of finding
\begin{equation}
\min_{W\in \mathbb R^2} \sum_{j=1}^n m_j |P_jW| \, ,
\label{GFTn}
\end{equation}
in the case of existence of solution $ W_{\ast}=(x_{\ast},y_{\ast}) \not\in \{P_j\}_{j=1}^n $, has the value (\ref{Wcost}) equal to
$$
\left|\sum_{j=1}^n m_jz_j(\cos \beta_j+ \mathbf i \sin \beta_j) \right|\stackrel{(\ref{WcostSum})}{=}
\left|\sum_{j=1}^n m_j(z_j-z_{\ast})(\cos \beta_j+ \mathbf i \sin \beta_j)\right|, \quad z_{\ast}:=x_{\ast} + \mathbf i y_{\ast} \, .
$$
If the direction $ (\cos \beta_j, \sin \beta_j ) $ of the terminal $ P_j $ is counted starting from the $ x $-axis, then $ 2\pi -\beta_j $ is just the argument of the number $ z_j-z_{\ast} $.

However, the posed restriction on the network topology is significant for the constructive computation of the directions of the terminals in the network. Indeed, in the case of networks possessing the claimed property, formulas for the direction of any terminal can be represented by radicals via the weights of the problem provided that these directions are counted starting from any edge of the optimal network. Note that the coordinates of facilities are not required for this aim, the only assurance of the existence of the network in the prescribed topology matters.

This is not the case for the networks with facilities of the degree higher than $ 3 $. Even the problem (\ref{GFTn}) cannot be resolved by radicals for
general case of $ n=4 $ terminals. For instance, when the weight $ m $ from Example \ref{ExWeber4} increases continuously from the initial value $ m=4 $, the facilities $ W_1 $ and $ W_2 $ of the corresponding optimal bifacility networks tend to a collision point $ W_{\ast} $ at the value $ m=m_{\ast} $.
The point $ W_{\ast} $ is the solution to the unifacility Weber problem (\ref{GFTn}) with $ n=4 $. Coordinates of $ W_{\ast} $ satisfy the $ 10 $th degree algebraic equations over $ \mathbb Z $; so do the critical value $ m_{\ast} \approx 4.326092 $ \cite{Uteshev_Semenova_2020_Ar}.

\section{Conclusions}

\textbf{1.} Though the deduced formulas for the length or cost of optimal networks do not require the coordinates of any network facility, it is possible to
generate them in the above exploited complexification ideology. For instance, in the complex plane, the Steiner point of the triangle
$ P_1P_2P_3 $ (in the case of its existence and provided that the triangle vertices are numbered counterclockwise) is given as
$$
\frac{\mathbf i}{\sqrt{3}} \left(\varepsilon_2z_2-\varepsilon_1z_1+ (\overline{z_1}-\overline{z_2})\mathcal L/\overline{\mathcal L} \right)
$$
where $ \mathcal L:= z_3+\varepsilon_2 z_1 + \varepsilon_1 z_2 $ is the expression from formula (\ref{Length3}), i.e. $ |\mathcal L| $ is the length of the full Steiner tree. For the terminals of Example \ref{Ex0}, this formula yields
$$
\frac{454+250\sqrt{3} + \mathbf i(262+ 150\sqrt{3})}{112+60\sqrt{3}}=\frac{731+95\sqrt{3}}{1703} +\mathbf i \frac{293+135\sqrt{3}}{1703}
\approx 4.108004+ \mathbf i \,  2.416637 \, .
$$

\vspace{1em}

\textbf{2.} It looks challenging to find  the minimum value for
$$
\left| \sum_{j=1}^n z_j \widetilde{U}_j \right|
$$
in the set of all vectors $ (\widetilde{U}_1,\dots, \widetilde{U}_n) \in \mathbb C^n $ such that
$$
\left\{ \widetilde{U}_j^6=1 \right\}_{j=1}^n \ , \ \sum_{j=1}^n  \widetilde{U_j}=0 \, .
$$
It resembles a knapsack problem.

\vspace{1em}

\setcounter{equation}{0}
\setcounter{theorem}{0}
\setcounter{example}{0}

\end{document}